\documentclass[a4paper]{amsart}
\usepackage{graphicx} % Required for inserting images

\usepackage{amsthm}

\usepackage{tikz-cd}
\usepackage{amscd}
\usepackage{hyperref}
\usepackage{mathrsfs}
\usepackage{xcolor}
\usepackage{setspace}
\usepackage{csquotes}
\usepackage[only, llbracket, rrbracket]{stmaryrd}
\usepackage{amssymb, amsmath, amsthm, mathrsfs, amsfonts, latexsym, bbm}

\usepackage{comment}

\usepackage[shortlabels]{enumitem}

\newtheorem{theorem}{Theorem}[section]
\newtheorem{corollary}[theorem]{Corollary}

\newtheorem{proposition}[theorem]{Proposition}

\newcounter{claim}
\newcommand{\norm}[1]{\left\lVert#1\right\rVert}
\newcommand{\abs}[1]{\left\lvert#1\right\rvert}

\title{Complemented ideals of $\ell_\infty$}
\author{Michael Hru\v s\' ak and Luis S\'aenz }
\date{August 2025}

\thanks{
{The authors gratefully acknowledge support from a DGAPA-PAPIIT grant IN101323 and  a \linebreak 
SECIHTI grant CBF-2025-I-898.}}

\address{Centro de Ciencas Matem\'aticas\\
UNAM\\
A.P. 61-3, Xangari, Morelia, Michoac\'an\\
58089, M\'exico}

\email{michael@matmor.unam.mx, luisdavidr@ciencias.unam.mx}

\begin{document}

\begin{abstract} Answering questions raised in \cite{Leonetti, Uzcategui} we characterize ideals $\mathcal I\subseteq \mathcal P(\omega)$ such that $c_{0,\mathcal I}$ is complemented in $\ell_\infty$ as exactly those ideals for which the space $K_{\mathcal I}= \mathsf{Stone}(\mathcal P(\omega)/\mathcal I)$ is approximable, i.e., the unit ball of the space $M(K_{\mathcal I})$ of  signed Radon measures on $K_\mathcal I$ is separable in the weak* topology.
\end{abstract}

\maketitle

The content of these notes concerns the structure of the ideals of $\ell_\infty$ - the Banach lattice of all bounded sequences of reals equipped with the supremum norm and pointwise order.  It is important to distinguish between the two kinds of ideals we deal with. 

Recall that a \emph{Banach lattice} is a Banach space equipped with a lattice order relation such that if $\abs{x}=x\lor-x$, then for any $x,y\in X$, if $\abs{x}\leq \abs{y}$ then $\norm{x}\leq \norm{y}$. An \emph{ideal} $Y$ of a Banach lattice $X$ is a linear subspace of $X$, such that if $y\in Y$, $x\in X$ and $\abs{x}\leq \abs{y}$, then $x\in Y$. 

Meanwhile, given a nonempty set $X$, $\mathcal I\subseteq \mathcal P (X)$ is an \emph{ideal} if $A,B\in \mathcal I$ implies that $A\cup B\in \mathcal I$, and if $A\in \mathcal I$ and $B\subseteq A$, then $B\in \mathcal I$. An ideal is \emph{proper} if $X\notin \mathcal I$. All ideals we consider consist of subsets of $\omega$ (the set of all natural numbers) and contain all finite subsets of $\omega$.

Closed Banach ideals of the Banach lattice $\ell_\infty$ and set theoretic ideals of $\mathcal P(\omega)$ are closely related as we now explain. Given an ideal $\mathcal I\subseteq \mathcal P (\omega)$ and a sequence $(x_n)_{n\in \omega}\in \ell_\infty$ we say its $\mathcal I$-limit is $r\in \mathbb R$ if for any $\epsilon>0$, $\{n\in \omega: \abs{x_n-r}\ge \epsilon\}\in \mathcal I$, and denote this fact by $\mathcal I \text{-}\lim x_n=r$. So we may define 
$$c_{0,\mathcal I}=\{(x_n)_{n\in \omega}\in \ell_\infty: \mathcal I\text{-}\lim x_n=0\}.$$

As noted by Uzcategui and Rincón-Villamizar \cite{Uzcategui1}, every closed ideal (in the Banach lattice sense) $X\subseteq \ell_\infty$ that contains $c_0$  is of the form $c_{0,\mathcal I}$ for some (set theoretic) ideal $\mathcal I\subseteq \mathcal P (\omega)$.

There has been interest in the question: \emph{When is $c_{0,\mathcal I}$ complemented in $\ell_\infty$?} The history of this question started at the problem session of the 45th Winter School in Abstract Analysis (Czech Republic, 2017). Tomasso Russo communicated this question to Paolo Leonetti \cite{Leonetti}, who showed that if $\mathcal P (\omega)/\mathcal I$ is not ccc, then the answer is negative. Not long after Tomasz Kania \cite{Kalton} noticed that a little more can be proved: if $\mathcal P (\omega)/\mathcal I$ is not ccc, then $C(K_\mathcal I)$ cannot be embedded into $\ell_\infty$ (we clarify the notation in the following section). As a consequence of our main result, both answers are equivalent, i.e., $c_{0,\mathcal I}$ is complemented in $\ell_\infty$ if and only if $C(K_{\mathcal I})$ can be embedded into $\ell_\infty$. Recently, Carlos Uzcategui and Michael A. Rincón-Villamizar \cite{Uzcategui} proposed a series of properties to analyze the question. In the last section, we provide the pertinent (counter-)examples.

\medskip

Our main result seems to give a definitive answer to these questions:

\begin{theorem}\label{main}
    Let $\mathcal I\subseteq\mathcal P (\omega)$ be an ideal. Then the following are equivalent:
    \begin{enumerate}
        \item  $c_{0,\mathcal I}$ is complemented in $\ell_\infty$.
        \item $K_{\mathcal I}$ is approximable.
        \item $C(K_\mathcal I)$ is isometric to a subspace of $\ell_\infty$.
        \item $C(K_{\mathcal I})$ is isomorphic to a subspace of $\ell_\infty$.
        \item  There is an $\epsilon\in (0,1)$ such that $(\mathcal P (\omega)/\mathcal I)^+$ can be covered by a countable family of subsets each of them with intersection number bigger than $1-\epsilon$.
    \end{enumerate}
\end{theorem}

In particular, ideals satisfying the conditions of the theorem have to be such that $\mathcal P(\omega)/\mathcal I$ admits a strictly positive finitely additive measure, while  ideals such that $\mathcal P(\omega)/\mathcal I$ is $\sigma$-centered, i.e. $K_\mathcal I$ is separable, provide examples.

\section{Preliminaries}

The following are standard notions in their respective areas.  

\smallskip

Given an ideal $\mathcal I$, we denote, as usual, by $\mathcal I ^*$ its dual filter; i.e., $\mathcal I^*=\{\omega\setminus I:I\in \mathcal I\}$. We denote by $K_\mathcal I=\{\mathcal U\in \beta\omega:\mathcal I^*\subseteq \mathcal U\}$ the corresponding compact subspace of $\beta\omega$. 

We say that $\mathbb{B}$   is \emph{ccc} if every antichain in $\mathbb{B}$ is countable and it is \emph{$\sigma$-centered} if $\mathbb{B}^+$ can be covered by countably many ultrafilters.
 By $\mathsf{Stone}(\mathbb B)$ we denote the Stone space of the algebra $\mathbb B$ (i.e the unique compact space such that $\mathbb B\simeq \mathsf{Clopen(K)}$, the Boolean algebra of clopen sets.) A compact space $K$ is \emph{extremally disconnected} if the closure of any open set is open, equivalently if the Boolean algebra of clopen subsets of $K$ is complete.

\smallskip
We say that a Boolean algebra $\mathbb{B}$   \emph{admits a strictly positive measure} if there is a finitely additive measure $\mu$ on $\mathbb{B}$  such that for any $b\in\mathbb{B}^+=\mathbb{B}\setminus \{\bf 0\}$, $\mu(b)>0$. Given $s\in \mathbb B^{<\omega}$ we denote $i(s)=\max\{\abs{F}:F\subseteq \text{dom}(s)\land \cap_{i\in F}s(i)\neq \varnothing\}$. Given $\mathcal A\subseteq \mathbb B$, the \emph{intersection number} of $\mathcal A$ is defined as:

$$I(\mathcal A)=\inf\Bigg\{\frac{i(s)}{\abs s}:s\in\mathcal A^{<\omega}\Bigg\}.$$

\smallskip

Recall Kelley's fundamental result:

\begin{theorem}[Kelley \cite{Kelley}]\label{Kelley}
    Let $\mathbb{B}$ be a Boolean algebra. $\mathbb{B}$ admits a strictly positive measure if and only if $\mathbb{B}^+$ can be covered by a countable family of subsets with positive intersection number.
\end{theorem}

Given a compact space $K$, a \emph{Radon measure} over $K$ is a regular finite measure over the Borel algebra of $K$. We denote by $M(K)$ the set of all Radon signed measures over a compact space $K$. This set, equipped with the total variation norm is a Banach space. By Riesz representation theorem $M(K)$ can be naturally regarded as $C(K)^{*}$, the dual Banach space of the space of all continuous functions over $K$. Following Plebanek and D\v{z}amonja \cite{Mirna}, we say $K$ is  \emph{approximable} if the unit ball of $M(K)$ (the set of all measures of norm at most $1$) is separable in the $w^\ast$-topology. This property may be characterized strictly in terms of the Boolean algebra:

\begin{theorem}[Mägerl-Namioka \cite{Mager}]\label{Mager-Namioka} Let $\mathbb B$ be a Boolean algebra. The space $\mathsf{Stone}(\mathbb B)$ is approximable if and only if there is $\epsilon\in (0,1)$ such that $\mathbb{B}^+$ can be covered by a countable family of subsets each of them with intersection number bigger than $1-\epsilon$.
    
\end{theorem}

We call a Banach space $X$ \emph{injective} if for any  Banach spaces $Y\subseteq Z$  and any operator $T:Y\to X$, there exists an operator $\hat T:Z\to X$ that extends $T$ (i.e. $\hat T\lvert_Y=T$). Furthermore, we call $X$ \emph{$1$-injective} if $\hat T$ can be found such that $\frac{\norm{\hat T}}{\norm{T}}=1$.

\begin{theorem}[Goodner-Nachbin, see \cite{Kalton}]
Let $K$ be a compact space.
    $C(K)$ is $1$-injective if and only if $K$ is extremally disconnected.
\end{theorem}

\begin{proposition}[\cite{Kalton}]
    Let $X$ be an injective Banach space, and let $Y$ any Banach space such that $X$ embeds into $Y$. Then $X$ is isomorphic to a complemented subspace of $Y$.
\end{proposition}

Given two Banach spaces $X$ and $Y$, recall that $X$ embeds into $Y$ if there is a linear injective operator with closed range $T:X\to Y$, and that $X$ embeds isometrically into $Y$ if there is a linear operator $T:X\to Y$ such that $\norm{T(x)}=\norm x$. In this context, we will use the following result by Talagrand \cite[Theorem II.1]{Talagrand}:

\begin{theorem}[\cite{Talagrand}]\label{Talagrand}
    Let $K$ be a compact space. $C(K)$ embeds into $\ell_\infty$ if and only if $C(K)$ embeds isometrically into $\ell_\infty$.
\end{theorem}

Finally, recall the fundamental result of Lindenstrauss:

\begin{theorem}[Lindenstrauss, see \cite{Kalton}] \label{Kalton}
    Let $X$ be an infinite dimensional Banach subspace of $\ell_\infty$. Then $X$ is complemented if and only if $X\simeq \ell_\infty$.
\end{theorem}

\section{The main results}

In order to prove Theorem \ref{main}, we require the following folklore result. We provide its short proof.

\begin{proposition}\label{complete BA}
    Let $\mathcal I\subseteq \mathcal P (\omega)$ be an ideal such that $\mathcal P (\omega)/\mathcal I$ is ccc, then $\mathcal P (\omega)/\mathcal I$ is a complete Boolean algebra.
\end{proposition}
\begin{proof}
    It is enough to prove that if $\mathcal A,\mathcal B\subseteq \mathcal P (\omega)/\mathcal I$ are disjoint such that $\mathcal A\cup\mathcal B$ is a maximal antichain, then they can be \emph{separated}, i.e., there is a $c\in\mathbb B$ such that $a\leq c$ for every $a\in \mathcal A$ while $b\land c =\bf 0$ for every $b\in \mathcal B$. 
    
    By ccc, both families are countable, as there are no $(\omega,\omega)$-gaps in $\mathcal P(\omega)/[\omega]^{<\omega}$, therefore there is $C\in [\omega]^{\omega}$ that separates $(\mathcal A,\mathcal B)$.
    
    If $\mathcal A=\varnothing$, then $C$ separates. Otherwise, if $A\in \mathcal A$, then  $C\in\mathcal I ^+$; because $A \in \mathcal I ^+$ and $A\subseteq^* C$.
\end{proof}

The next result was originally announced in \cite{Talagrand} without a proof. We provide its easy proof.

\begin{proposition}\label{Approximable}
    Let $K$ be a compact space. $C(K)$ is isometric to a subspace of $\ell_\infty$ if and only if $K$ is approximable.
\end{proposition}
\begin{proof}
  Assume that there is an isometry $T:C(K)\to\ell_\infty$. Let $(x^*_n)_{n\in \omega}$ be any $w^*$-dense subset of  $B_{\ell_\infty^*}$. Consider $\mu_n$ as the (signed) measures defined by $\mu_n=T^*(x^*_n)$. As $T$ is an isometry, $T^*$ is also an isometry. Therefore, for every $n$, $\norm{\mu_n}\le 1$. The $w^*$-density of $(x^*_n)_{n\in \omega}$ implies that $(\mu_n)_{n\in \omega}$ is $w^*$-dense in $B_{M(K)}$.

  Assume that $K$ is approximable. Let $(\mu_n)_{n\in \omega}$ be a $w^*$-dense sequence in $B_{M(K)}$ (i.e. the set of measures with norm at most 1), and let $T:C(K)\to \ell_\infty$ be defined by $T(f)(n)=\int_{K}fd\mu_n$. Linearity is clear.  Notice that
  $$\norm{f}=\sup_{n\in \omega}\abs{\mu_n(f)}=\sup_{n\in \omega}\abs{\int_{K}fd\mu_n}=\norm{T(f)}.$$

\end{proof}

We also need a result mentioned by Kania without a proof.

\begin{proposition}[Kania \cite{Kania}]\label{Kania}
    Let $\mathcal I\subseteq \mathcal P (\omega)$ be an ideal, then $\ell_\infty/c_{0,\mathcal I}$ is isometric to $ C(K_\mathcal I)$.
\end{proposition}
\begin{proof}
Consider $T:\ell_\infty\to C(\beta\omega)$ given by 
$T((x_n)_{n\in \omega})(\mathcal U)=\mathcal U\text{-}\lim x_n$. This is clearly a linear isometry. 
Notice $f\in T[c_{0,\mathcal I}]$ if and only if, for any $\mathcal U\in K_\mathcal I$ $f(\mathcal U)=0$. Indeed,
    $T((x_n)_{n\in \omega})\in T[c_{0,\mathcal I}]$ if and only if, $\mathcal I\text{-}\lim x_n=0$, if and only if, for any $\mathcal U\in K_{\mathcal I}$, $\mathcal U\text{-}\lim x_n=0$. This establishes that if $U:\ell_\infty/c_{0,\mathcal I}\to C(K_\mathcal I)$ is defined by $U((x_n)_{n\in \omega})+c_{0,\mathcal I})(\mathcal U)=\mathcal U\text{-}\lim x_n$, then the map is well defined and linear.

    Now consider the canonical quotient map $\pi: C(\beta\omega)\to C(\beta\omega)/T[c_{0,\mathcal I}]$. It is clear that $\norm{\pi(f)}=\sup\{\abs{f(\mathcal U)}:\mathcal U \in K_\mathcal I\}.$ As $T$ is an isometry, this establishes that $U$ is an isometry. Finally, an appeal to Tietze's extension theorem shows $U$ is surjective.
\end{proof}

\begin{proof}[Proof of Theorem \ref{main}]
($2\leftrightarrow 3 \leftrightarrow 4 \leftrightarrow 5$) Follows from Proposition \ref{Approximable}, Theorem \ref{Mager-Namioka}, and Theorem \ref{Talagrand}.
   
    ($1\rightarrow 3$)
    Assume that $c_{0,\mathcal I}$ is complemented in $\ell_\infty$. Then there is a closed subspace $Y\subseteq\ell_\infty$ such that $Y\oplus c_{0,\mathcal I}=\ell_\infty$. By Proposition \ref{Kania} we get that $C(K_\mathcal I)$ is isometric to $ Y $.

   ($2\to 1$) If $K_{\mathcal I}$ is approximable then by Theorem \ref{Mager-Namioka} and Theorem \ref{Kelley}, $\mathcal  I^+$ admits a strictly positive. Therefore, $\mathcal P(\omega)/\mathcal I$ is ccc. By Proposition \ref{complete BA} it is extremally disconnected. By Proposition \ref{Approximable} $C(K_\mathcal I)$ is isometric to a subspace of $\ell_\infty$. As $K_{\mathcal I}$ is extremally disconnected, $C(K_\mathcal I)$ is $1$-injective and so it is embedded as a complemented subspace, therefore, $c_{0,\mathcal I}$ is complemented.
\end{proof}

\textbf{Remark.} It is worth noting that for compact spaces of the form $K_{\mathcal I}$, $C(K_\mathcal I)$ being complemented in $\ell_\infty$ is equivalent with $C(K_\mathcal I)$ being (isometrically) isomorphic to a subspace of $\ell_\infty$, a fact that is not true for arbitrary compact spaces.
As a direct corollary of Theorem \ref{main} and Theorem \ref{Kalton} we get:

\begin{corollary}
    Let $\mathcal I\subseteq\mathcal P (\omega)$ be an ideal such that $C(K_\mathcal I)$ is infinite dimensional. Then we can add to the previous theorem the equivalent condition that $C(K_\mathcal I)\simeq\ell_\infty$.
\end{corollary}

Next we shall provide illuminating examples that together with Theorem \ref{main} answer several questions from \cite{Uzcategui}. 

The following observation appears in \cite[Corollary 1.2.]{Dow} attributed to Alan Dow:

\begin{proposition}[\cite{Dow}]\label{Dow}
    Let $\mathbb B$ be a complete Boolean algebra of size at most $\mathfrak c$, then there is an ideal $\mathcal I\subseteq \mathcal P(\omega)$ such that $\mathbb B\simeq \mathcal P(\omega)/\mathcal I$.
\end{proposition}

 Theorem \ref{Mager-Namioka} provides a combinatorial characterization
  of the fact that $C(\mathsf{Stone}(B))$ embeds into $\ell_\infty$.
  
\begin{proposition}[\cite{Plebanek}]
    Let $\mathbb{B}$ be a Boolean algebra such that $P(\mathsf{Stone}(\mathbb B))$ is approximable, then $\mathbb B$ admits a strictly positive measure.
\end{proposition}

As the standard measure algebra $\mathbb{B}(\omega_1)$ over $2^{\omega_1}$ to add $\omega_1$-many random reals does not fulfill the hypothesis of Theorem \ref{Mager-Namioka}, this exemplifies why the converse implication does not follow.

Recall that the Maharam type of a measure is the minimal size of a subset of the $\sigma$-algebra of measurable sets dense in the metric $d(A,B)=\mu(A\Delta B)$.

\begin{proposition}[\cite{Plebanek}]
    Let $K$ be a compact space and $\mu$ be a positive Radon measure over $K$ of countable Maharam type, then $K$ is approximable.
\end{proposition} 

The authors of \cite{Mirna} provide a {\sf ZFC} example of a Boolean algebra whose Stone space is approximable, but admits no strictly positive measure with countable Maharam type (a result that improves an earlier construction of \cite{Talagrand} under {\sf CH}).

\medskip

\textbf{Example 1.} \textit{There is an ideal $\mathcal I$ such that $\mathcal P(\omega)/\mathcal I$ is ccc but $C(K_\mathcal I)$ does not embed into $\ell_\infty$.}

For instance, a classical example in \cite{Gaifman} presents a Boolean algebra $\mathbb B'$ that is ccc, but does not admit a strictly positive measure. Let $\mathbb B$ be its completion. Find, using Proposition \ref{Dow}, an ideal such that $\mathcal P (\omega)/\mathcal I\simeq \mathbb B$. Then $K_\mathcal I$ cannot be approximable, therefore $C(K_\mathcal I)$ cannot be embedded into $\ell_\infty$.

\medskip

\textbf{Example 2.} \textit{There is an ideal $\mathcal I$ such that $c_{0,\mathcal I}$ is complemented but  $\mathcal P(\omega)/\mathcal I$ is not $\sigma$-centered ($\omega$-maximal in the language of \cite{Uzcategui}).}

\smallskip

Consider the \emph{random forcing} $\mathbb B(\omega)={\mathsf {Borel}(2^\omega)}/\mathsf{null}$. This is a complete Boolean algebra whose space of probability measures is $w^*$-separable (as the Maharam type of the Lebesgue measure is countable) but is not $\sigma$-centered. Use \ref{Dow} once again to find the desired example.

\medskip

We want to point out that our main result sheds some light on Problem $9$ of \cite{OnGrothendick}: \emph{Characterize those ideals $\mathcal I$ for which $c_{0,\mathcal I}$ is a Grothendieck space\footnote{Recall that a Banach space $X$ is \emph{Grothendieck} if every sequence of elements of $X^*$ that is $w^*$-convergent is $w$-convergent.}}. It is known that any complemented subspace of a Grothendieck space is itself Grothendieck. As $\ell_\infty$ is Grothendieck, any ideal $\mathcal I$ that fulfills any (and therefore all) of the properties of our main theorem is such that $c_{0,\mathcal I}$ is Grothendieck.

Finally, we would like to know the answers to the following questions:
\begin{enumerate}[label=\alph*)]
    \item \textit{
Is it possible to construct a ZFC example of an extremally disconnected compact space $K$ such that it is approximable but does not admit a strictly positive measure of countable Maharam type?}

\item \textit{Is it possible to give a ZFC example of an extremally disconnected compact space $K$ such that the ball of $M(K)$ is $w^*$-separable but $M(K)$ is not?} \footnote{Counterexamples to both properties but without the condition of being extremally disconnected can be found in \cite{Mirna} and \cite{AvilesPlebanekRodriguez}, respectively.}
\end{enumerate}

We find these questions interesting because \textit{negative} answers would allow us to add more properties to the list of equivalent conditions of our main theorem: (a) there is a measure on $K$ of countable Maharam type, or (b) $M(K)$ is $w^*$-separable. 

\bigskip

{\bf Acknowledgment.} We would like to thank the anonymous referee for making several suggestions that greatly improved the quality of the paper. 

%\printbibliography

\bibliographystyle{plain}
\bibliography{references}

\end{document}